\documentclass[a4paper,10pt]{article}

\usepackage[utf8]{inputenc}
\usepackage[english]{babel} 

\usepackage[all,cmtip]{xy}
\usepackage{setspace}
\usepackage{multirow,rotating,graphicx}
\DeclareGraphicsExtensions{.jpg}

\usepackage{amsmath,amsfonts,amssymb}
\usepackage{amsthm}
\usepackage{booktabs}  
\usepackage{wasysym}
\usepackage{array}
\usepackage{xtab}
\usepackage{amscd} 
\usepackage{latexsym,mathtools}
\usepackage{ytableau}  
\usepackage[a4paper]{geometry}
\usepackage{enumitem}

\usepackage{color}
\usepackage{rotating}
\usepackage{tikz,fancyvrb} 
\usetikzlibrary{arrows.meta,arrows,automata}
\usepackage{cancel}
\usepackage{datetime}
\usepackage{subfig}
\usepackage{adjustbox}
\usepackage{caption}  
\usepackage[section]{placeins}

\usepackage[big]{layaureo} 
\usepackage[colorlinks, citecolor=blue]{hyperref}
\usepackage{frcursive}

\theoremstyle{plain}

\newtheorem{te}{Theorem}[section]

\newtheorem{theorem}[te]{Theorem}
\newtheorem{lemma}[te]{Lemma}

\newtheorem{proposition}[te]{Proposition}
\theoremstyle{definition}

\newtheorem{definition}[te]{Definition}

\theoremstyle{remark}

\newtheorem{question}[te]{Question}

\newtheorem{remark}[te]{Remark}

\newcommand{\field}[1]{\mathbb{#1}} %
\newcommand{\R}{\field{R}} 

\newcommand{\lie}[1]{{\mathfrak{#1}}}     
\newcommand{\g}{{\lie{g}}}

\newcommand{\SO}{\mathrm{SO}}

\newcommand{\so}{\mathfrak{so}}

\newcommand{\hook}{\lrcorner\,}  

\DeclareMathOperator{\dd}{\textsl{d}}
\DeclareMathOperator{\tr}{tr}        
\DeclareMathOperator{\Id}{Id}        
\DeclareMathOperator{\Span}{span}
\DeclareMathOperator{\Ric}{Ric}      
\DeclareMathOperator{\ric}{ric}      
\DeclareMathOperator{\ad}{ad}        
\DeclareMathOperator{\Tr}{tr}        
\DeclareMathOperator{\Der}{Der}      

\numberwithin{equation}{section}

\allowdisplaybreaks

\newcolumntype{C}{>{$}c<{$}}
\newcolumntype{L}{>{$}l<{$}}
\newcolumntype{R}{>{$}r<{$}}

\title{A non-Standard Indefinite Einstein Solvmanifold}

\author{Federico A. Rossi}

 \date{\today\ \currenttime} 

\begin{document}
\maketitle

\begin{abstract}
We describe an example of an indefinite invariant Einstein metric on a solvmanifold which is not standard, and whose restriction on the nilradical is nondegenerate.
\end{abstract}
%
\setcounter{tocdepth}{2}
\tableofcontents

\renewcommand{\thefootnote}{\fnsymbol{footnote}}
\footnotetext{\emph{MSC class 2020}: \emph{Primary} 53C50; \emph{Secondary} 53C25, 53C30, 22E25}
\renewcommand{\thefootnote}{\arabic{footnote}}

\section{Introduction}

The study of Einstein metrics on homogeneous spaces has been a central topic in differential geometry, with deep implications for both mathematics and physics. A particular area of focus within this field is the study of solvmanifolds, namely connected, simply connected solvable Lie groups equipped with left-invariant metrics. These manifolds offer a fertile ground for exploring different geometric structures, including Einstein metrics, i.e. metrics such that the Ricci curvature is proportional to the metric itself:
\begin{equation}
\Ric_g=\lambda\Id,\quad \lambda\in\R.
\end{equation}
While extensive research has been conducted on the Riemannian case, where the metric is positive definite, the indefinite case -- characterized by metrics that admit both positive and negative eigenvalues -- presents additional complexities that have only recently begun to be systematically explored.

In the Riemannian setting, the classification and construction of Einstein metrics on solvmanifolds have been extensively investigated. A foundational result, due to the work of J.~Heber, J.~Lauret and Y.~Nikolayevsky \cite{Heber:noncompact,Lauret:Einstein_solvmanifolds,Nikolayevsky}, is the correspondence between Einstein metrics and nilsolitons, which is a special class of Ricci solitons defined on nilpotent Lie groups. The subsequent work of C.~B\"ohm and R.~Lafuente~\cite{Bohm_Lafuente_2023}, J.~Lauret~\cite{Lauret:RicciSoliton} and C.~Will~\cite{LauretWill:EinsteinSolvmanifolds}, M.~Jablonski \cite{Jablonski1,Jablonski2} has further expanded this theory, providing a comprehensive understanding of the algebraic and geometric structures underlying these spaces. In particular, an Einstein solvmanifold $\tilde{\g}$ is standard, meaning that it decomposes as $\tilde{\g}=[\tilde{\g},\tilde{\g}]\oplus^{\perp}\lie{a}$, where $\lie{a}$ is an abelian subalgebra of $\tilde{\g}$. Moreover, $\tilde{\g}$ is of Iwasawa-type, roughly speaking, for any $X\in\lie{a}$, $\ad X$ is symmetric with respect to the metric. Excellent references on the positive-definite setting include the survey by J.~Lauret \cite{Lauret:SurveyNilsoltions} and the more recent one by  M.~Jablonski \cite{Jablonski:Survey}.

However, this correspondence does not directly extend to the indefinite case and the situation changes dramatically when considering indefinite metrics, which introduce significant complexities (for example, the presence of null directions). In this context, the existence and construction of Einstein metrics are far less understood, demanding innovative approaches and techniques. Together with D.~Conti, we first studied the Einstein metrics on nilpotent Lie algebras \cite{ContiRossi:EinsteinNilpotent,ContiRossi:EinsteinNice,ContiRossi:RicciFlat,ContiDelBarcoRossi:DiagramInvolutions}, and later we focused on the relations between Einstein solvmanifolds and indefinite nilsolitons \cite{ContiRossi:IndefiniteNilsolitons,ContiRossi:NiceNilsolitons} (see also the works of Z.~Yan \cite{Yan:Pseudo-RiemannianEinsteinhomogeneous} and S.~Deng \cite{YanDeng:DoubleExtensionRiemNilsolitons}). These results are related to the so called pseudo-Iwasawa Einstein solvmanifolds, which are the indefinite Einstein metrics most similar to the Riemannian ones. More recently, in joint works with R.~Segnan Dalmasso, we explored some different Einstein metrics: in~\cite{ContiSegnanRossi:PseudoSasaki,ContiSegnanRossi:ErratumPseudoSasaki,ContiSegnanRossi:PseudoSasakiEinsteinSolv} we studied pseudo-Sasaki Einstein solvmanifolds,  and in~\cite{ContiSegnanRossi:NotBased} we constructed Einstein metrics independent of nilsolitons. Both of these constructions are examples of non-pseudo-Iwasawa Einstein solvmanifolds. These works underline the fundamental disparities between the Riemannian and indefinite cases. The particularly curious reader should find a brief and partial description of the indefinite Einstein solvmanifolds in~\cite{Rossi:NewSpecialPseudoEinstein}.

In the pioneering paper~\cite{ContiRossi:IndefiniteNilsolitons} it was claimed that there exists an Einstein solvmanifold $\tilde{\g}$ admitting an orthogonal decomposition $\tilde{\g}=[\tilde{\g},\tilde{\g}]\oplus^{\perp}\lie{h}$ such that $[\tilde{\g},\tilde{\g}]$ is  the nilradical and $\lie{h}$ is nonabelian (\cite[Example~1.7]{ContiRossi:IndefiniteNilsolitons}). However, this example turned out to be incorrect, as the metric is not Einstein. This leaves an open question: whether it is possible to find such an example in the indefinite setting. More precisely we ask  the following:
\begin{question}\label{question}
Does there exist a solvmanifold $\tilde{\g}$ with an Einstein left-invariant metric $\tilde{g}$ admitting an orthogonal decomposition $\tilde{\g}=\g\oplus^{\perp}\lie{h}$, where $\g=[\tilde{\g},\tilde{\g}]$ is the nilradical, the restriction of the metric on $\g$ is nondegenerate, and $\lie{h}$ is nonabelian?
\end{question}
In the present article, we answer this problem affirmatively (Theorem~\ref{thm:Example}), correcting the mistake of our previous work.

This paper aims to explore and elucidate the fundamental differences between the Riemannian and indefinite cases of Einstein solvmanifolds. By synthesizing recent advancements in this field, we seek to provide a clearer understanding of the unique challenges posed by indefinite Einstein metrics. Through this comparative study, we hope to contribute to the ongoing development of pseudo-Riemannian geometry and to the broader understanding of Einstein metrics in more general settings.

More precisely, in Section~\ref{sec:preliminaries} we recall the basic definitions and general setting. In the next Section~\ref{sec:IndefiniteVSDefinit} we briefly recap the differences between indefinite and definite Einstein solvmanifolds, and in the last Section~\ref{sec:Example} we describe an example of a solvmanifold $\tilde{\g}$ with an Einstein metric $\tilde{g}$ that does not admit a standard decomposition (Theorem~\ref{thm:Example}): in fact, we can decompose the Lie algebra as $\tilde{\g}=[\tilde{\g},\tilde{\g}]\oplus^{\perp}\lie{h}$, but the subspace $\lie{h}$ is not abelian. Moreover, we also observe that classical Azencott-Wilson results cannot be applied to find a pseudo-Iwasawa solvmanifold isometric to $\tilde{\g}$ (Remark~\ref{rem:AzencottWilsonNonFunziona}), and that none of the elements in $\lie{h}$ can be used to satisfy a general nilsoliton equation in the sense of~\cite{ContiSegnanRossi:NotBased} (Remark~\ref{rem:NotGenNil}), emphasizing the uniqueness and peculiarity of this example.

\bigskip
\noindent \textbf{Acknowledgments:}  
The author would like to thank the Organizing Committee and the Scientific Committee of the International Scientific Conference ``Algebraic and Geometric Methods of Analysis'', May 27-30, 2024, Ukraine.
The author acknowledge GNSAGA of INdAM and the PRIN project n. 2022MWPMAB ``Interactions between Geometric Structures and Function Theories''.

\section{Invariant Structures on Lie Groups}\label{sec:preliminaries}

We recall some definitions and concepts that we will use throughout the paper.

Given a Lie group $G$, we denote its Lie algebra by $\g$. We focus our attention on \emph{left-invariant pseudo-Riemannian metrics} $g$ defined on a Lie group $G$. These metrics can be represented as bilinear nondegenerate symmetric forms $g$ (i.e., definite or indefinite scalar products) on the corresponding Lie algebra $\g$. The pair $(\g,g)$ is referred to as a \emph{metric Lie algebra}. Consequently, the Levi-Civita connection $\nabla$, the Riemann curvature $R$, the Ricci operator $\Ric$, and the Ricci curvature $\ric$ are also left-invariant and can be expressed as linear tensors on $\g$. A metric is classified as \emph{Riemannian} if it is positive definite; otherwise, it is considered an \emph{indefinite metric}. The \emph{nilradical} of a Lie algebras is its the maximal nilpotent ideal.

The \emph{lower central series} of a Lie algebra $\g$ is defined recursively as $\g^0:=\g$, $\g^i:=[\g,\g^{i-1}]$. A Lie algebra is called \emph{nilpotent} if its lower central series converges to the trivial subspace $\{0\}$, indicating the existence of a natural number $s$ such that $\g^s=\{0\}$. In this case, the minimum such $s$ is known as the \emph{step} of the Lie algebra.

The \emph{derived series} of $\g$ is defined as $\mathfrak{a}_0:=\g$, $\mathfrak{a}_i:=[\mathfrak{a}_{i-1},\mathfrak{a}_{i-1}]$. If there exists a natural number $r$ such that $\mathfrak{a}_r=\{0\}$, the Lie algebra is said to be \emph{solvable}. The dimension of a complementary subspace of $\g'=[\g,\g]=\lie{a}_1$ is known as the \emph{rank}. For example, a rank-one solvable Lie algebra can be expressed as a direct sum of its derived subalgebra $\g'$ and the span of a single element $X\in\g$. 

A Lie algebra is considered \emph{unimodular} if the trace of the adjoint operator $\ad (X)$ is zero for all $X\in\g$. It is evident that a nilpotent Lie algebra is also solvable and unimodular.

A metric is classified as \emph{Einstein} if its Ricci tensor $\ric$ is proportional to the metric itself, i.e., $\ric= \lambda g$ for some constant $\lambda\in\R$. Equivalently, the Ricci operator $\Ric$ is a scalar multiple of the identity. If the proportionality constant $\lambda$ is zero, the metric is said to be \emph{Ricci-flat}. The \emph{scalar curvature} is defined as the trace of the Ricci operator and is denoted by $s:=\tr(\Ric)$. For Einstein metrics with nonzero scalar curvature, we have the relation $s=n\lambda\neq0$.

A \emph{solvmanifold} $(G,g)$ is defined as a simply connected solvable Lie group $G$ equipped with a left-invariant metric $g$. 


To describe a Lie algebra, we will use the notation introduced by S.~Salamon \cite{Salamon:ComplexStructures}, which is now becoming popular and widely accepted. 
For example, the $4$-dimensional $2$-step nilpotent Lie algebra $\g$ will be described by giving the action of the Chevalley-Eilenberg operator $\dd$ on the dual (which is equivalent to giving the expression for the Lie bracket):
\[(0,0,e^{12},e^{13}),\]
which means that  $\g^*$ has a fixed basis $\{e^1,\dotsc, e^4\}$ with $\dd e^1=\dd e^2=0$, $\dd e^3=e^{12}=e^1\wedge e^2$ and $\dd e^4=e^{13}=e^1\wedge e^3$. This implies that $\{e_1,e_2,e_3,e_4\}$ is  a basis of the Lie algebra  and the nontrivial brackets are
\[[e_1,e_2]=-e_3,\qquad [e_1,e_3]=-e_4.\]

The following result (see~\cite{ContiRossi:EinsteinNilpotent}) is useful for computing the Ricci tensor of a metric Lie algebra $\g$:
\begin{lemma}[{\cite[Lemma 1.1]{ContiRossi:EinsteinNilpotent}}]
The Ricci tensor $\ric$ of a left-invariant pseudo-Riemannian metric~$g$ on a Lie algebra $\g$ is given by:
\begin{multline}
\label{eqn:ricciingeneral}
\ric(v,w)=\frac12g  (\dd v^\flat,  \dd w^\flat) -\frac12g (\ad v, \ad w)\\
-\frac12 \Tr\ad(v\hook  \dd w^\flat+w\hook  \dd v^\flat)^\sharp-\frac12 B(v,w),\quad v,w\in\g
\end{multline}
where $B(v,w)=\tr(\ad v\circ\ad w)$ is the Killing form.
\end{lemma}

In the Riemannian setting, Einstein solvmanifolds are closely related to nilsoliton metrics on their nilradical (see e.g., \cite{Heber:noncompact,Lauret:Einstein_solvmanifolds,Lauret:SurveyNilsoltions} or the Section~\ref{sec:IndefiniteVSDefinit}), i.e., a metric $g$ on a nilpotent Lie algebra $\g$ that satisfies 
\begin{equation}\label{eqn:nilsoliton}
\Ric=\lambda\Id +D
\end{equation}
for some derivation $D$ of $\g$. In the indefinite  setting, there is a more general condition that the restriction of an Einstein metric may satisfy on a nilpotent ideal. We pose the following definition (see~\cite{ContiSegnanRossi:PseudoSasakiEinsteinSolv})
\begin{definition}
On a nilpotent Lie algebra $\g$, a metric $g$ and a derivation $D$ satisfy the \emph{generalized nilsoliton equation} if
\begin{equation}
\label{eqn:generalizednilsolitoningeneral}
\Ric=\tau\biggl(- \Tr ((D^s)^2) \Id- \frac12[D,D^*] +(\Tr D)D^s\biggr), \quad \Tr (\ad v\circ D^*)=0,\ v\in\g,
\end{equation}
for some $\tau=\pm1$. Here $D^s=\frac12(D+D^*)$ is the symmetric part of $D$.
\end{definition}
Metrics that are generalized nilsoliton may be extended to rank-one Einstein  solvmanifolds (see \cite[Proposition~2.1]{ContiSegnanRossi:PseudoSasakiEinsteinSolv}). Observe that if $D=D^*$, then we recover the usual nilsoliton equation~\eqref{eqn:nilsoliton}.

Finally, we report a result stating a sufficient condition to have an isometry between two Lie algebras of the form $\g\ltimes\lie{a}$. This result belongs to the ``Azencott-Wilson tricks'' to prove isometry between Lie groups, since the first result appear in~\cite{AzencottWilson}. The following statement is the one in \cite[Proposition~1.1]{ContiSegnanRossi:PseudoSasakiEinsteinSolv} which is more general (but see also \cite[Theorem~5.6]{AzencottWilson} and \cite[Section~1.8]{EberleinHeber} for the Riemannian case, \cite[Proposition~1.19]{ContiRossi:IndefiniteNilsolitons} for an ``indefinite standard'' version or \cite[Proposition~2.2 and Corollary~2.3]{ContiSegnanRossi:PseudoSasaki} without the standard assumption on indefinite metric).

\begin{proposition}[{\cite[Proposition~1.1]{ContiSegnanRossi:PseudoSasakiEinsteinSolv}}]
\label{prop:pseudoAzencottWilson}
Let $H$ be a subgroup of $\SO(r,s)$ with Lie algebra $\lie h$ and $\widetilde{\g}$ a Lie algebra of the form $\widetilde{\g}=\g\rtimes \lie{a}$ endowed with a $H$-structure. Let $\chi\colon\lie a\to\Der(\g)$ be a Lie algebra homomorphism such that, extending $\chi(X)$ to $\widetilde{\g}$ by declaring it to be zero on $\lie a$,
\begin{equation}
\label{eqn:adXstar}
\chi(X)-\ad X\in\lie h, \qquad [\chi (X),\ad Y]=0,\ X,Y\in\lie a.
\end{equation}
Let $\widetilde{\g}^*$ be the Lie algebra $\g\rtimes_\chi\lie{a}$. If $\widetilde{G}$ and $\widetilde{G}^*$ denote the connected, simply connected Lie groups with Lie algebras $\widetilde{\g}$ and $\widetilde{\g}^*$, with the corresponding left-invariant $H$-structures, there is an isometry from $\widetilde{G}$ to $\widetilde{G}^*$, whose differential at $e$ is the identity of $\g\oplus\lie{a}$ as a vector space, mapping the $H$-structure on $\widetilde{G}$ into the $H$-structure on $\widetilde{G}^*$.
\end{proposition}

\begin{remark}
Throughout this paper, we will primarily deal with left-invariant objects on the Lie group $G$, identifying them with their corresponding linear counterparts on the Lie algebra $\g$. Our primary objective is to investigate Einstein metrics on solvable Lie groups. It is important to note that tensors defined on $\g$ naturally correspond to left-invariant tensors on the connected, simply-connected Lie group $G$ with Lie algebra $\g$, thereby allowing us to translate our subsequent discussions into the language of Lie groups.
\end{remark}

\section{Definite vs. Indefinite Einstein Solvmanifolds}\label{sec:IndefiniteVSDefinit}

In \emph{Riemannian geometry}, an important category of metric Lie algebras is comprised of standard solvable Lie algebras. These are solvable Lie algebras $\tilde\g$ equipped with a fixed metric such that the orthogonal complement of their derived subalgebra, $[\tilde \g,\tilde \g]^\perp$, is abelian. Explicitly, \emph{standard Riemannian} Lie algebras can be decomposed as
\begin{equation}
\label{eqn:standardRiemannian}
\tilde \g=\g\oplus^\perp\lie{a},
\end{equation}
where $\g= [\tilde \g,\tilde \g]$ is nilpotent and $\lie{a}$ is abelian. The symbol $\oplus^\perp$ denotes an orthogonal direct sum of vector spaces, and the abelian subalgebra $\lie{a}$ acts nontrivially unless $\tilde \g$ is itself abelian.

It was known that Riemannian Einstein solvmanifolds are nonunimodular~\cite{Dotti:RicciCurvature}, and a classical result by J.~Milnor~\cite{Milnor:curvatures} states that the only nilpotent Lie algebra with an Einstein left-invariant metric is the abelian one, which is Ricci-flat.

In fact, J.~Lauret~\cite{Lauret:Einstein_solvmanifolds} established that all Riemannian Einstein solvmanifolds satisfy~\eqref{eqn:standardRiemannian}, i.e., they are standard. Furthermore, the structural characteristics of standard Riemannian Einstein solvmanifolds had been previously elucidated by J.~Heber in~\cite{Heber:noncompact}. Combining  the work of J.~Heber and J.~Lauret together, we can conclude that in an Einstein solvmanifold $\tilde{\g}$ the nilradical coincides with the derived algebra $[\tilde{\g},\tilde{\g}]$ (see~\cite{Lauret:SurveyNilsoltions}).

Indeed, J.~Heber established that the abelian subalgebra $\lie a$ operates on $\g$ via normal derivations. More specifically, for any element $X\in\lie a$, the metric adjoint of $\ad X$ is a derivation that commutes with $\ad \lie a$.

Subsequent research by R.~Azencott and E.N.~Wilson~(\cite{AzencottWilson}) demonstrates that the Lie algebra can be modified through projection onto its self-adjoint component, resulting in an isometric solvmanifold where the abelian subalgebra $\lie a$ acts by self-adjoint derivations. Furthermore, it has been shown in~\cite{Heber:noncompact} that the action of $\lie a$ is faithful and that there exists an element $H$ in $\lie a$ such that the adjoint operator $\ad H$ is positive definite. Consequently, the Lie algebra $\g$ is characterized as being of \emph{Iwasawa type}. More precisely, the \emph{Iwasawa type} conditions as stated in \cite{Heber:noncompact,Wolter:EinsteinSolvable} are the following:
\begin{enumerate}[label=(Iw\arabic*)]
\item\label{Iwa1} The orthogonal complement $\lie{a}$ of $[\tilde{\g},\tilde{\g}]$ is abelian.
\item\label{Iwa2} For all  $X\in \lie{a}$, $\ad X$ are symmetric relative to the metric $\tilde  g$ and nonzero for $X\neq0$.
\item\label{Iwa3} For some $H\in\lie{a}$, the restriction $\ad H$ on $[\tilde{\g},\tilde{\g}]$ is positive definite.
\end{enumerate}

In  fact, it is known (\cite{Heber:noncompact}) that the restriction of the metric on $[\tilde{\g},\tilde{\g}]$ is a \emph{nilsoliton metric}, meaning that its Ricci operator satisfies the equation $\Ric=\lambda\Id+D$ for some derivation $D$ of $[\tilde{\g},\tilde{\g}]$). The element $H$ from the previous Condition~\ref{Iwa3} is (a multiple) of the derivation $D$ arising from the nilsoliton condition. This derivation $D$ is referred to as the \emph{pre-Einstein} or \emph{Nikolayevsky derivation} of $[\tilde{\g},\tilde{\g}]$, and it was proven in~\cite{Nikolayevsky} that $D$ is diagonalizable and satisfies the property $\tr(D\circ X)=\tr(X)$ for any derivation $X$ of $[\tilde{\g},\tilde{\g}]$.

\smallskip

However, the situation is quite different when the Einstein metric is not positive definite. In~\cite{ContiRossi:IndefiniteNilsolitons}, motivated by the Riemannian cases, we introduced the following  definitions for  indefinite metrics:

\begin{definition}
\label{def:standard}
A \emph{standard decomposition} of a pseudo-Riemannian metric Lie algebra $\tilde \g$ is a decomposition
\[\tilde \g=\g\oplus^\perp\lie a,\]
where $\g$ is a nilpotent ideal and $\lie a$ is an abelian subalgebra.
\end{definition}

\begin{definition}
Given a solvable Lie algebra $\tilde \g$, we say that a standard decomposition $\tilde\g=\g\oplus^\perp\lie{a}$ is \emph{pseudo-Iwasawa} if
\begin{equation}
 \label{eqn:pseudo-Iwasawa}
 \ad X= (\ad X)^*, \quad X\in\lie{a}.
\end{equation}
\end{definition}

These definitions generalized the usual  Riemannian  ones: note that standard Riemannian Lie algebras admits a standard decomposition, and that Riemannian Lie algebras  of Iwasawa  type are pseudo-Iwasawa.

We will briefly recall some differences between the indefinite and the definite setting, with  particular attention to the structure of the algebras. A detailed analysis of the differences in the nilsoliton condition or the metric can be found in~\cite{ContiRossi:IndefiniteNilsolitons} or~\cite{ContiSegnanRossi:NotBased}, and a short summary on results about the indefinite Einstein solvmanifolds can be found in~\cite{Rossi:NewSpecialPseudoEinstein}. 

Assume $\tilde{\g}$ is an Einstein solvmanifold with an indefinite metric $\tilde{g}$. Then we can make the following remarks:
\begin{enumerate}[label=(D\arabic*)]
\item There exist Einstein metrics on nonabelian nilpotent Lie algebras: while Ricci-flat examples were known in the literature, the first non-Ricci-flat examples were found in~\cite{ContiRossi:EinsteinNilpotent},  and more examples and constructions were obtained in~\cite{ContiRossi:RicciFlat,ContiRossi:EinsteinNice,ContiDelBarcoRossi:DiagramInvolutions,Conti:RicciFlatness}.
\item Examples of Einstein metrics on nilpotent Lie algebras with scalar curvature $s\neq0$ are unimodular, in contrast with the Riemannian case (see~\cite{Milnor:curvatures} and~\cite{Dotti:RicciCurvature}).
\item While existence of Einstein metrics with nonzero scalar curvature on a nilpotent Lie algebra is obstructed by an algebraic condition \cite[Theorem~4.1]{ContiRossi:EinsteinNilpotent}, Ricci-flat metrics on nilpotent Lie algebras exist in great abundance~\cite{Conti:RicciFlatness}.
\item Taking Einstein metric on a nilpotent Lie algebra $\g$, we see that the derived algebra $[\g,\g]$ does not necessarily coincide with the nilradical (see,  for example, \cite[Example~1.2]{ContiRossi:IndefiniteNilsolitons} or examples constructed in~\cite{ContiRossi:RicciFlat,ContiRossi:EinsteinNice,ContiDelBarcoRossi:DiagramInvolutions}).
\item The restriction of the metric may be degenerate on the nilradical (see e.g., \cite[Example~1.3]{ContiRossi:IndefiniteNilsolitons}), or it may be  degenerate on the derived algebra (\cite[Example~1.1]{ContiRossi:IndefiniteNilsolitons}).
\item\label{DifferenceNonStandard} There are Einstein solvmanifolds that do not admit any standard decomposition as in Definition~\ref{def:standard}; see \cite[Example~1.6]{ContiRossi:IndefiniteNilsolitons} where the nilradical coincides with the derived algebra, and the restriction of the metric is degenerate.
\item The notions of algebraic Ricci soliton and nilsoliton can be naturally translated to the pseudo-Riemannian setting. It is known (see~\cite{Onda:ExampleAgebraicRicciSolitons}) that pseudo-Riemannian algebraic Ricci solitons are Ricci solitons; but there are Lie groups carrying left-invariant indefinite metrics which are Ricci solitons without satisfying~\eqref{eqn:nilsoliton}, see~\cite{BatatOnda:AlgebraicRicciSoliton}.
\item In indefinite signature, there exists $4$ different types of nilsolitons, i.e. metric such that $\Ric=\lambda\Id+D$ for some derivation $D$. The derivation $D$ may not be related to the Nikolayevsky derivation, and may not be diagonalizable (see \cite[Section~2]{ContiRossi:IndefiniteNilsolitons}).
\item Different nilsoliton metrics may appear on the same nilpotent Lie algebra (see \cite[Remark~2.6]{ContiRossi:IndefiniteNilsolitons} or \cite{KondoTamaru:LorentzianNilpotent}). This shows that the uniqueness of nilsoliton metrics up to scaling and automorphisms as proved in~\cite{Lauret:RicciSoliton} does not hold for indefinite metrics.
\item In indefinite signature, a more general nilsoliton equation~\eqref{eqn:generalizednilsolitoningeneral} may be used to construct Einstein solvmanifold (\cite[Proposition~2.1]{ContiSegnanRossi:PseudoSasakiEinsteinSolv}, \cite{ContiSegnanRossi:NotBased}), and the restriction of the metric on a nilpotent ideal does not necessarily define a nilsoliton equation as in the Riemannian case.
\item\label{DifferenceNonPseudoIwasawa} There are Einstein solvmanifolds with a standard decomposition that are not of pseudo-Iwasawa type. Some examples can be found in \cite[Example~1.14, Example~3.2 or Example~3.12]{ContiRossi:IndefiniteNilsolitons}, but a more systematic construction of Einstein solvmanifolds non-pseudo-Iwasawa is given by the pseudo-Sasaki Einstein solvmanifolds presented in~\cite{ContiSegnanRossi:PseudoSasakiEinsteinSolv} (in fact, by \cite[Proposition~2.6]{ContiSegnanRossi:PseudoSasaki}, a Sasaki pseudo-Riemannian solvmanifold does not admit any pseudo-Iwasawa decomposition).
\item Observation~\ref{DifferenceNonPseudoIwasawa} shows that Condition~\ref{Iwa2} may fail for indefinite metrics.
\item In~\cite{ContiRossi:NiceNilsolitons} we extend nilsoliton metrics to rank-one Einstein solvmanifolds, even if the eigenvalues of the Nikolayevsky derivation are nonpositive. Consequently, Condition~\ref{Iwa3} is not satisfied.
\item\label{DifferenceLast} Observation~\ref{DifferenceNonStandard} demonstrates that Condition~\ref{Iwa1} fails in the indefinite setting if we consider an Einstein solvmanifold with a degenerate metric on the nilradical.
\end{enumerate}

The last observation~\ref{DifferenceLast} indicates that there exist  non-standard Einstein solvmanifolds; however, in the known examples, the restriction of the metric is degenerate on the nilpotent ideal $\g$.

It remains to investigate the existence problem described in Question~\ref{question}, namely whether it is possible to construct an Einstein solvmanifold $\tilde{\g}$ with a decomposition $\tilde{\g}=\g\oplus^{\perp}\lie{h}$ where $\lie{h}$ is not abelian and the  metric on $\g$ is nondegenerate. The existence of such an example was first claimed in \cite[Example~1.7]{ContiRossi:IndefiniteNilsolitons}, but that example turned out to be incorrect because the metric was not Einstein. Nevertheless, the claim remains valid, as we will describe in detail in the next Section~\ref{sec:Example}.

\section{A non-Standard Indefinite Einstein Solvmanifold}\label{sec:Example}

In this  section,  we answer positively to Question~\ref{question} with an explicit example.

Consider the  $5$-dimensional Lie algebra $\tilde{\g}$ with structure equation given by
\begin{equation}\label{LieAlgebraEx}
\left( e^{42}+e^{51}-e^{54} , -e^{41}+e^{52} , e^{12}-e^{51}+2e^{53}-\frac{7}{12}e^{54} , 0 , 0 \right).
\end{equation}
Explicitly, the bracket in the basis $\{e_1,\dots e_5\}$ are the following:
\begin{equation}\label{eq:brkRules}
\begin{gathered}
[e_1,e_2]=-e_3,\qquad [e_4,e_2]=-e_1,\qquad [e_4,e_1]=e_2,\qquad  [e_5,e_1]=-e_1+e_3,\\
[e_5,e_2]=-e_2,\qquad [e_5,e_3]=-2e_3,\qquad [e_5,e_4]=e_1+\frac{7}{12}e_3.
\end{gathered}
\end{equation}

We start with the
\begin{lemma}
The nilradical $\lie{n}$ of $\tilde{\g}$ is given by $\lie{n}=\Span\{e_1,e_2,e_3\}$.
\end{lemma}
\begin{proof}
We note  that the derived algebra is 
\[[\tilde{\g},\tilde{\g}]=\Span\{e_1,e_2,e_3\},\]
and the nilradical satisfies $\lie{n}\supseteq[\tilde{\g},\tilde{\g}]$. However, the only linear combination $\alpha e_4+\beta e_5$ that give a nilpotent element  in the adjoint algebra is $0$, so we can conclude that the nilradical coincide with the derived algebra.
\end{proof}

Now we will consider the invariant indefinite metric $\tilde{g}$ on $\tilde{\g}$ defined  by
\begin{equation}\label{eq:metricExample}
\tilde{g}=
\frac{497}{576}e^1\otimes e^1+\frac{49}{192} e^2\otimes e^2+2e^3\otimes e^3-\frac{7}{6}e^1\odot e^3
-\frac{245}{6144}e^4\otimes e^4 -\frac{1225}{6144}e^5\otimes e^5 ,
\end{equation}
and let $\hat{g}$ be the restriction of $\tilde{g}$ on the nilradical $\lie{n}$. 
\begin{lemma}
The restriction $\hat{g}$ defines a Riemannian metric on $\lie{n}$, in particular $\hat{g}$ is not degenerate.
\end{lemma}
\begin{proof}
Choosing on $\g$ the orthogonal basis 
\[\hat{e}_1=e_1,\qquad \hat{e}_2=e_2,\qquad \hat{e}_3=96 e_1+71 e_3\]
one easily sees that $\hat{g}(\hat{e}_i,\hat{e}_i)>0$ for  all $i=1,2,3$.
\end{proof}
By extending the orthogonal basis of  the  previous Lemma with $\hat{e}_4=e_4$ and $\hat{e}_5=e_5$, we see that $\tilde{g}$ is an indefinite metric of signature $(3,2)$ on $\tilde{\g}$.

Let $\lie{h}$ be the orthogonal complement of the nilradical $\lie{n}$ on $\tilde{\g}$, i.e. $\lie{h}=\lie{n}^{\perp}$.
\begin{lemma}\label{lem:lieh}
It  holds $\lie{h}=\lie{n}^{\perp}=\Span\{e_4,e_5\}$.
\end{lemma}
\begin{proof}
Since $\Span\{\hat{e}_4,\hat{e}_5\}=\Span\{e_4,e_5\}$, the statement follows using the orthogonal basis $\{\hat{e}_1,\dots,\hat{e}_5\}$ described above.
\end{proof}

\begin{lemma}\label{lem:hnonabelian}
The space $\lie{h}$ is not contained in any abelian subalgebra of $\tilde{\g}$.
\end{lemma}
\begin{proof}
By Lemma \ref{lem:lieh} we have that $\{e_4,e_5\}$ is  a basis of $\lie{h}$ as vector space, and  by the bracket rules~\eqref{eq:brkRules} we  get
\[[e_5,e_4]=e_1+\frac{7}{12}e_3\neq0.\qedhere\]
\end{proof}

\begin{lemma}\label{lem:NoStandardDecomposition}
There is no nilpotent ideal $\g\subseteq \tilde{\g}$ and no abelian subalgebra  $\lie{a}\subseteq\tilde{\g}$ such that the following orthogonal decomposition holds
\[\tilde{\g}=\g\oplus^{\perp}\lie{a}.\]
\end{lemma}
\begin{proof}
Assume by contradiction that there exists a nilpotent ideal $\g$ of $\tilde{\g}$ such that $\tilde{\g}=\g\oplus^{\perp}\lie{a}$ and $\lie{a}$ is abelian.

Since the nilradical  $\lie{n}$ is the maximal nilpotent ideal of $\tilde{\g}$ and the metric $\hat{g}$ on $\lie{n}$ is positive definite, it follows that $\g\subseteq\lie{n}$ and the restriction of the metric $\tilde{g}$ on $\g$ is positive definite and nondegenerate. Hence, using Lemmas~\ref{lem:lieh} we obtain that \[\lie{a}=\g^{\perp}\supseteq\lie{n}^{\perp}=\lie{h}=\Span\{e_4,e_5\}\]
and by Lemma~\ref{lem:hnonabelian} we conclude that $\lie{a}$ is not abelian, which is an absurd.
\end{proof}

A direct computation gives:
\begin{lemma}\label{lem:EinsteinMetricExample}
The Ricci tensor $\widetilde{\Ric}$ of the metric $\tilde{g}$ is 
\[\widetilde{\Ric}=\frac{4096}{175}\sum_{i=1}^5 e^i\otimes e^i=\frac{4096}{175}\Id,\]
i.e., $\tilde{g}$ is an invariant Einstein metric on $\tilde{\g}$ with  Einstein constant  $\lambda=\frac{4096}{175}$.
\end{lemma}
\begin{proof}
It follows directly from a straightforward computation using~\eqref{eqn:ricciingeneral}.
\end{proof}

Using the previous Lemmas, we can state the following properties:
\begin{theorem}\label{thm:Example}
The invariant metric $\tilde{g}$ defined in~\eqref{eq:metricExample} 
on the solvable Lie algebra $\tilde{\g}$ with structure equation~\eqref{eq:brkRules}
is an Einstein metric with Einstein constant $\lambda=\frac{4096}{175}\neq0$.

Moreover,  the metric Lie algebra $\tilde{\g}$ does not  admit any standard decomposition.
\end{theorem}
\begin{proof}
The first part follows from Lemma~\ref{lem:EinsteinMetricExample}.

The second part is the statement of  Lemma~\ref{lem:NoStandardDecomposition}.
\end{proof}

It  is natural  to ask whether it is possible to find a pseudo-Iwasawa Einstein solvmanifold that is isometric to the solvable Lie algebra $\tilde{\g}$ described in~\eqref{LieAlgebraEx}. The following Remark~\ref{rem:AzencottWilsonNonFunziona} shows that it is not possible to apply Proposition~\ref{prop:pseudoAzencottWilson} for this purpose.

\begin{remark}\label{rem:AzencottWilsonNonFunziona}
As shown earlier, $\lie{h}$ is generated by $\{e_4,e_5\}$. Now, assume by contradiction that there exists a map $\chi\colon\lie a\to\Der(\g)$ as in Proposition~\ref{prop:pseudoAzencottWilson}, and let $X=\chi(e_5)$ be the derivation of $\g$ such that $[X,\ad e_4]=0$. Easy computations show that a generic derivation of $\g$ can be expressed, with respect to the basis $\{e_1,e_2,e_3\}$, as the matrix
\begin{equation}\label{genericoDerivationg}
\begin{pmatrix}
 a_{1,1} & a_{1,2} & 0 \\
 a_{2,1} & a_{2,2} & 0 \\
 a_{3,1} & a_{3,2} & a_{1,1}+a_{2,2} \\
\end{pmatrix},\qquad a_{i,j}\in\R.
\end{equation}
With  respect to the same basis, we have $\ad e_4=\begin{pmatrix}
 0 & -1 & 0 \\
 1 & 0 & 0 \\
 0 & 0 & 0 \\
\end{pmatrix}$. Thus, if $X$ commutes with $\ad e_4$, $X$ must have the form:
\[\chi(e_5)=X=\begin{pmatrix}
 x_{1,1} & x_{1,2} & 0 \\
 -x_{1,2} & x_{1,1} & 0 \\
 0 & 0 & 2 x_{1,1} \\
\end{pmatrix},\qquad x_{i,j}\in\R.\]
However, since $\chi(e_5)-\ad  e_5$ must belong to $\lie{h}\subset\so(3,0)$, it follows that $\chi(e_5)-\ad  e_5$ must be antisymmetric with respect to the metric $g$. This condition can be written as:
\[\chi(e_5)-\ad  e_5=-\chi(e_5)^*+(\ad  e_5)^*\iff \chi(e_5)+\chi(e_5)^*=\ad  e_5+(\ad  e_5)^*.\]
A straightforward computation yields to
\begin{gather*}
\ad e_5 = 
\begin{pmatrix}
 -1 & 0 & 0 \\
 0 & -1 & 0 \\
 1 & 0 & -2 \\
\end{pmatrix},
\qquad 
(\ad  e_5)^*=
\begin{pmatrix}
 -\frac{14}{3} & 0 & \frac{32}{7} \\
 0 & -2 & 0 \\
 \frac{1}{36} & 0 & -\frac{4}{3} \\
\end{pmatrix},\\
X^*=
\begin{pmatrix}
	-\frac{41}{15} x_{1,1}   & -\frac{7}{5} x_{1,2}   & \frac{32}{5} x_{1,1}  \\[2pt]
	\frac{71}{21} x_{1,2}    & x_{1,1}                & -\frac{32}{7} x_{1,2} \\[2pt]
	-\frac{497}{180} x_{1,1} & -\frac{49}{60} x_{1,2} & \frac{86}{15} x_{1,1},
\end{pmatrix},\\
X+X^*=
\begin{pmatrix}
	-\frac{26}{15} x_{1,1}   & -\frac{2}{5} x_{1,2}   & \frac{32}{5}x_{1,1}    \\[2pt]
	\frac{50}{21}x_{1,2}     & 2 x_{1,1}              & -\frac{32}{7}  x_{1,2} \\[2pt]
	-\frac{497}{180} x_{1,1} & -\frac{49}{60} x_{1,2} & \frac{116}{15} x_{1,1}
\end{pmatrix}.
\end{gather*}
This lead to no solution for the equation  $X+X^*=\ad  e_5+(\ad  e_5)^*$, which  is  a contradiction.
\end{remark}

Finally, we observe that no element in $\lie{h}$ satisfies the generalized nilsoliton equation~\eqref{eqn:generalizednilsolitoningeneral} for the  restriction $g$ of the Einstein metric on $\tilde{\g}$:

\begin{remark}\label{rem:NotGenNil}
A generic element in $\lie{h}=\Span{e_4,e_5}$ acts on $\g$ as a linear combination of $\ad e_4$  and $\ad_5$. Thus,  we can  consider the derivation $D$ of $\g$ written, in the basis $\{e_1,e_2,e_3\}$, as
\[D=\alpha \ad e_4+\beta \ad e_4=
\begin{pmatrix}
	-\beta & -\alpha & 0        \\
	\alpha & -\beta  & 0        \\
	\beta  & 0       & -2 \beta
\end{pmatrix},\quad \alpha,\beta\in\R.\]
A direct computation shows that $-\tr ((D^s)^2) \Id - \frac{1}{2} [D, D^*] + (\tr D)  D^s$ equals
\[
\begin{pmatrix}
	\frac{14 \beta ^2}{3}-\frac{50 \alpha ^2}{21}  & \frac{8}{15} \alpha  \beta                    & \frac{16}{7} \left(\alpha ^2-5 \beta ^2\right)        \\[2pt]
	\frac{176}{21} \alpha \beta                    & -\frac{2 \alpha ^2}{5}-\frac{62 \beta ^2}{21} & -\frac{96}{7}  \alpha  \beta                          \\[2pt]
	\frac{49 \alpha ^2}{120}+\frac{13 \beta ^2}{8} & -\frac{259}{180} \alpha  \beta                & -\frac{2}{105} \left(73 \alpha ^2+345 \beta ^2\right)
\end{pmatrix}.
\]
On the other hand, using~\eqref{eqn:ricciingeneral}, we find that the Ricci  operator of $g$ is
\[\Ric=
\begin{pmatrix}
	-\frac{36864}{1715} & 0                   & 0                  \\
	0                   & -\frac{36864}{1715} & 0                  \\
	-\frac{6144}{245}   & 0                   & \frac{36864}{1715}
\end{pmatrix}.\]
Hence, there are no solutions  to \eqref{eqn:generalizednilsolitoningeneral}, either with $\tau=+1$ or $\tau=-1$.
\end{remark}

\bibliographystyle{plainurl}
\bibliography{A_nonStandard_ExampleBIB}

\small\noindent Dipartimento di Matematica e Informatica, Universit\`a degli Studi di Perugia, via Vanvitelli 1, 06123 Perugia, Italy.\\
\texttt{federicoalberto.rossi@unipg.it}

\newenvironment{dedication}
  {\clearpage           
   \thispagestyle{empty}
   \vspace*{\stretch{1}}
   \itshape             
   \raggedleft          
  }
  {\par 
   \vspace{\stretch{3}} 
   \clearpage           
  }

\begin{dedication}
{\large \cursive{Dedicato alla mia Mamma,\\
 al suo infinto amore,\\
 al suo enorme affetto.}}
\end{dedication}

\end{document}